\documentclass[11pt,leqno]{article}

\usepackage{amsmath,amsfonts,amscd,amssymb,theorem}

\long\def\comment#1\endcomment{}

\makeatletter
\begingroup
\gdef\th@dotted{\normalfont\itshape
  \def\@begintheorem##1##2{%
        \item[\hskip\labelsep \theorem@headerfont ##1\ ##2.]}%
\def\@opargbegintheorem##1##2##3{%
   \item[\hskip\labelsep \theorem@headerfont ##1\ ##2\ (##3).]}}
\endgroup
\makeatother

\theoremstyle{dotted}

\newtheorem{theorem}{Theorem}[section]
\newtheorem{lemma}[theorem]{Lemma}

\makeatletter
\begingroup
\gdef\th@upshape{\normalfont
  \def\@begintheorem##1##2{%
        \item[\hskip\labelsep \theorem@headerfont ##1\ ##2.]}%
\def\@opargbegintheorem##1##2##3{%
   \item[\hskip\labelsep \theorem@headerfont ##1\ ##2\ (##3).]}}
\endgroup
\makeatother

\theoremstyle{upshape}

\newtheorem{defn}[theorem]{Definition}
\newtheorem{remark}[theorem]{Remark}
\newtheorem{exa}[theorem]{Example}

\makeatletter
\renewcommand{\subsection}{\@startsection{subsection}{2}{0pt}{-3ex
plus -1ex minus -0.2ex}{-2mm plus -0pt minus
-2pt}{\normalfont\bfseries}} \makeatother

\makeatletter
\@addtoreset{equation}{section}
\makeatother

\newcommand{\cntrct}                
{\hspace{2pt}\raisebox{1pt}{\text{$\lrcorner$}}\hspace{2pt}}

\newcommand{\proof}[1][Proof.]{\smallskip\noindent{\em #1}}
\def\endproof{\hfill\ensuremath{\square}\par\medskip}

\def\eqref#1{\thetag{\ref{#1}}}

\let\latexref=\ref
\def\ref#1{{\normalfont{\latexref{#1}}}}

\newcommand{\wt}{\widetilde}

\newcommand{\idot}{{\:\raisebox{1pt}{\text{\circle*{1.5}}}}}
\newcommand{\hdot}{{\:\raisebox{3pt}{\text{\circle*{1.5}}}}}
\newcommand{\calo}{{\cal O}}
\newcommand{\Spec}{\operatorname{Spec}}

\newcommand{\Fun}{\operatorname{{\mathcal F}{\it un}}}

\newcommand{\sspec}{\operatorname{{\cal S}{\it pec}}}

\newcommand{\gr}{\operatorname{\sf gr}}

\newcommand{\C}{{\mathcal C}}
\newcommand{\A}{{\mathcal A}}

\newcommand{\T}{{\mathcal T}}
\newcommand{\N}{{\mathcal N}}

\newcommand{\Ext}{\operatorname{Ext}}

\newcommand{\rhom}{\operatorname{{\cal R}{\cal H}{\it om}}}
\newcommand{\RHom}{\operatorname{RHom}}

\newcommand{\ch}{\operatorname{ch}}

\newcommand{\Coker}{\operatorname{{\sf Coker}}}

\newcommand{\E}{{\cal E}}
\newcommand{\F}{{\cal F}}
\newcommand{\J}{{\cal J}}

\newcommand{\id}{\operatorname{\sf id}}
\newcommand{\Id}{\operatorname{\sf Id}}
\newcommand{\ev}{\operatorname{\sf ev}}
\newcommand{\tr}{\operatorname{\sf tr}}

\newcommand{\Tor}{\operatorname{Tor}}
\newcommand{\Shv}{\operatorname{Shv}}

\newcommand{\m}{{\mathfrak m}}

\title{Hochschild homology and Gabber's Theorem}

\author{D. Kaledin\thanks{Partially supported by CRDF grant
RUM1-2694-MO05.}}

\begin{document}

\maketitle

\section*{Introduction}

About twenty-five O. Gabber proved his famous theorem \cite{G} which
claims, roughly speaking, that the singular support $SS(\F)$ of a
$D$-module $\F$ on a smooth algebraic manifold $M$ is an {\em
involutive} subvariety in the cotangent bundle $T^*M$. Involutive
here means that the ideal sheaf $\J$ defining $SS(\F)$ is closed
with respect to the natural Poisson bracket on $T^*M$ -- that is,
$\{\J,\J\} \subset \J$.

This statement would have been very easy and straightforward, but
there a complication: $SS(\F)$ here is taken with reduced scheme
structure.

To understand the diffuculty, let us recall the definition of
$SS(\F)$. One first finds a good increasing filtration on $\F$ which
in particular is compatible with the order filtration on the algebra
$D$ of differential operators. The associated graded quotient $\gr
\F$ then becomes a module over $\gr D$, and can be localized to a
sheaf on $T^*M = \sspec \gr D$. By definition, $SS(\F)$ is the
support of the sheaf $\gr \F$.  The annihilator $\J' \subset
\calo_{T^*M}$ is easily seen to be involutive; however, $\J'$ need
not be a prime ideal sheaf, and the actual ideal sheaf $\J \subset
\calo_{T^*M}$ defining $SS(\F) \subset T^*M$ is the radical of
$\J'$.

Thus the statement one has to prove is entirely algebraic, but
highly non-trivial: there is no obvious reason why the radical of
this particular involutive ideal should also be involutive (and it
is certainly not true for all involutive ideals, e.g. the square of
{\em any} ideal is involutive for trivial reasons). The general
algebraic conjecture on involutivity of the singular support first
appeared in \cite{GS}. In the years between \cite{GS} and \cite{G},
a lot of progress was made -- following partial results in \cite{GS}
itself and \cite{bo}, a complete proof was given in \cite{KKS}. But
all these results used difficult analytic methods, such as
microlocalization and pseudodifferential operators of infinite
order. A different proof can be extracted from \cite{KS}, especially
Chapter 11, but it also uses difficult analytics facts. The argument
Gabber found was very beautiful and quite general, and it was purely
algebraic -- in fact, completely elementary in the sense that all
the techniques used are contained in an undergraduate algebra
course. However, the conceptual essence of the argument seems to be
very hard to catch. Perhaps as the result of this, the standard
textbooks on $D$-module theory such as \cite{borel} avoid it by
using a trick of J. Bernstein (although some explanation is given in
\cite{GC}).

Recently considerable progress has been made in some areas of
homological algebra, and we believe that one can now revisit
Gabber's Theorem and put it into a conceptual framework. From our
point of view, the appropriate notions are those of Hochschild
Homology and Cohomology for a small abelian category, and the
associated deformation theory.

Unfortunately, at present these theories are under development, and
far from being completed. The situation with Hochschild homology is
better, since a very thorough treatment has been given by B. Keller
\cite{kel} (nevertheless, there are fine points here as well, see
Remark~\ref{tr.rem}). The notion of Hochschild cohomology is in fact
much simpler; however, developing a deformation theory for abelian
categories based on Hochschild cohomology is a highly non-trivial
matter, and it has been done only very recently by W. Lowen and
M. Van den Bergh \cite{LB1}, \cite{LB2}, \cite{L}. At the moment --
possibly because of my lack of competence -- it is not clear to me
whether the theory constructed in these papers is strong enough. We
note that the requirements for the Gabber's Theorem are not very
strong, but they are quite specific; when a general theory is being
developed, they are likely to be omitted at first.

The present paper arose as an attempt to explain and generalize
Gabber's Theorem by first developing the Hochschild cohomology
theory in an alternative way, which works in lesser generality than
\cite{kel} and \cite{LB1} but is better adapted to this particular
problem. However, it soon became clear that a complete treatment
would require some space. Therefore we have decided to split the
text. While the theoretical paper \cite{K1} is being prepared, the
present paper is intended to give the actual practical proof of the
Gabber's Theorem by using Hochschild homology, and to state clearly
what general facts about Hochschild homology and cohomology one
needs. Necessarily, to do this we have to quote some things without
proof. Nevertheless, we have decided that doing this might be
useful. If nothing else, this will show what statements should
definitely be contained in a Hochschild homology package suitable
for practical applications, and allow a reader who is prepared to
accept on faith some general nonsense to understand how the proof of
Gabber's Theorem really works.

The paper is organized as follows. In Section~\ref{gen} we describe
the general properties that one expects from the Hochschild homology
and cohomology formalism -- specifically, those properties that we
need for the Gabber's Theorem. So as not to leave the exposition
completely without a foundation, in Section~\ref{shv} we sketch a
skeleton theory satisfying these requirements for sheaves on a
(regular affine) scheme with supports in a closed subscheme. We
leave two crucial compatibility results without proof (in a general
theory, these should be more-or-less automatic). Finally, in
Section~\ref{pf} we show how this formal theory allows one to prove
Gabber's Theorem, and we indicate possible generalizations.

\subsection*{Acknowledgments.} This paper owes its existence to
V. Ginzburg, my first math teacher, who convinced me that it would
be worth writing. I am very grateful to him for this particular
suggestion, and for his general continued interest in my work and
constant support. I have benefited from opportunities to discuss the
subject with R. Bezrukavnikov, A. Beilinson, F. Bogomolov,
A. Braverman, B. Feigin, D. Kazhdan, A. Kuznetsov, M. Lehn,
D. Tamarkin, B. Tsygan and M. Verbitsky, to whom I am also sincerely
grateful. S. Kleiman kindly indicated to me a mistaken reference in
the first version of this paper and pointed out the correct
reference \cite{GS}, and T. Monteiro Fernandes attracted my
attention to \cite{bo} and the relevant parts of \cite{KS}. I thank
P. Schapira for further explanations about the history of Gabber's
Theorem and the relations between different proofs.

\section{General requirements.}\label{gen}

Fix a field $k$, and let $\C$ be a small abelian $k$-linear
category. In our applications, $\C$ will have finite global
homological dimension. At least under this assumption, and possibly
always, one associates to $\C$ a pair of a homology and a cohomology
theory, called Hochschild homology and Hochschild cohomology, which
have some natural properties. Let us list those properties.

We start with Hochschild cohomology $HH^\hdot(\C)$ which is easier
to understand. First of all, it is a graded-commutative algebra over
$k$. Now, we note that the category $\Fun(\C,\C)$ of right-exact
functors from $\C$ to itself is also an abelian $k$-linear
category. The identity functor $\Id_\C$ is an object in
$\Fun(\C,\C)$. Then the following must hold.
\begin{itemize}
\item There exists a natural algebra map $\ev:HH^\hdot(\C) \to
  \Ext^\hdot(\Id_\C,\Id_\C)$, where the $\Ext$-groups are computed
  in $\Fun(\C,\C)$.
\end{itemize}
In most practical cases, this natural map will be an isomorphism,
and one can just take it as the definition of
$HH^\hdot(\C,\C)$. However, in some situations things might be more
complicated (for finiteness reasons). In any case, the map should
always exist. In particular, for any object $A \in \C$, we have by
restriction a natural {\em evaluation map}
$$
\ev_A:HH^\hdot(\C) \to \Ext_{\C}^\hdot(A,A).
$$
This map is compatible with the algeba structure on both sides.

The main interest in Hochschild cohomology is in its relation to
deformation theory. Fortunately, for Gabber's Theorem we only need
one-parameter deformations of order $1$, so that there is no need to
study obstructions (which would require, among other things, putting
a Lie algebra structure on $HH^\hdot(\C)$). We restrict our
attention to the following situation. Let $\C'$ be an abelian
$k$-linear category equipped with a functorial endomorphism
$h:\Id_{\C'} \to \Id_{\C'}$ such that $h^2=0$, and let $\C \subset
\C'$ be the full abelian subcategory of objects $A \in \C'$
annihilated by $h$. Denote by $\iota:\C \to \C'$ the embedding
functor. It is easy to see that it admits a right-adjoint functor
$\tau:\C' \to \C$ given by $\tau(A) = A/h = \Coker(h_A) \in \C
\subset \C'$.

\begin{defn}\label{cat.def}
The category $\C'$ is a {\em first-order deformation} of the
category $\C$ if the left-adjoint functor $\tau:\C' \to \C$ admits
derived functors $L^\hdot\tau:\C' \to \C$, and for any $k \geq 1$,
the composition $L^k\tau \circ \iota$ is isomorphic to $\Id_\C$. An
object $A \in \C'$ is said to be {\em flat} with respect to $\C
\subset \C'$ if $L^k\tau(A)=0$ for $k \geq 1$.
\end{defn}

\begin{exa}\label{exa}
If $A$ is a $k$-algebra, and $\C$ is the category of left
$A$-modules, then the category $\C'$ of left modules over $A[h]/h^2
= A \otimes_k k[h]/h^2$ is obviously a first-order deformation of
$\C$. The functor $L^\hdot\tau:\C' \to \C$ is equal to
$\Tor^\hdot_{A[h]/h^2}(A,-)$. This generalizes immediately to a
sheaf $\A$ of $k$-algebras on a topological space $X$, and the
category of sheaves of left $\A$-modules on $X$.
\end{exa}

\begin{remark}
The isomorphism $L^k\tau \circ \iota \cong \Id_\C$ is a flatness
condition on $\C'$. One can show that it is in fact sufficient to
require it for $k=1$. The same goes for flat objects.
\end{remark}

By definition, derived functors $L^k\tau$ of the right-exact functor
$\tau$ come equipped with an additional structure: they all fit into
a complex $L^\hdot\tau$. In particular, trunctating this complex to
$L^{\leq 1}\tau$, we obtain by Yoneda a canonical class $t_{\C'} \in
\Ext^2(\tau,L^1\tau)$. Composing this with the exact functor
$\iota:\C \to \C'$, we get a class $\iota(t_{\C'}) \in \Ext^2(\tau
\circ \iota,L^1\tau \circ \iota) = \Ext^2(\Id_\C,\Id_\C)$. It is
this class that should paramterize deformations.
\begin{itemize}
\item To any first-order deformation $\C'$ of the category $\C$, one
  associates a Hochschild cohomology class $\Theta_{\C'} \in
  HH^2(\C)$ such that 
\begin{equation}\label{def.comp}
\ev(\Theta_{\C'})=\iota(t_{\C'}).
\end{equation}
\end{itemize}

\begin{lemma}\label{no.ev}
Assume given a first-order deformation $\C' \supset \C$ and an
object $A \in \C$. If there exists a flat object $\wt{A} \in \C'$
such that $\wt{A}/h \cong A$, then $\ev_A(\Theta_{\C'})
=\iota(t_{\C'})_A =0$.
\end{lemma}

\proof{} To see explicitly the class $\iota(t_{\C'})_A$, we note
that for any short exact sequence $0 \to B \to C \to A \to 0$ in
$\C'$, we have a four-term exact sequence
\begin{equation}\label{4t}
\begin{CD}
0 @>>> \Coker p @>>> \tau(B) @>>> \tau(C) @>>> \tau(A)
@>>> 0,
\end{CD}
\end{equation}
where $p$ is the natural map $L^1\tau(C) \to L^1\tau(A)$. It is this
sequence that represents by Yoneda the class $p \circ t_{\C',A} \in
\Ext^2(\tau(A),\Coker p)$. Now, we apply this to the natural
exact sequence $0 \to N \to \wt{A} \to \iota(A) \to 0$, where $N$
the kernel of the surjective adjunction map $\wt{A} \to
\iota(A)$. Then since $\wt{A}$ is flat, $p$ is an isomorphism, and
the sequence \eqref{4t} becomes
$$
\begin{CD}
0 @>>> L^1\tau(\wt{A}) \cong A @>>> \tau(N) @>{0}>> A @>>> A @>>> 0,
\end{CD}
$$
which obviously represents $0$.
\endproof

In fact, the converse to this statement should also be true;
moreover, a full theory would show that up to an equivalence,
first-order deformations $\C'$ are classified by the corresponding
deformation classes $\Theta_{\C'} \in HH^2(\C)$. But we will not
need this.

Let us now turn to homology. Hochschild homology of the category
$\C$ should be a graded module $HH_\idot(\C)$ over the algebra
$HH^\hdot(\C)$. Fortunately, a very thorough theory of Hochschild
(and cyclic) homology for abelian categories has been developed by
B. Keller \cite{kel}; in particular, it is definitely known what are
the groups $HH_\idot(\C)$. Keller's theory even works in a more
general setting of exact categories, which allows to prove results
such as invariance with respect to derived equivalences and so
on. Unfortunately, Keller only works with homology, with no mention
of cohomology. Therefore the $HH^\hdot(\C)$-module structure on
$HH_\idot(\C)$ is not explicitly contained in his work. We will need
this structure, and will in fact need more.
\begin{itemize}
\item For any object $A \in \C$, there exists a natural trace map
$$
\tr_A:\Ext^\hdot_\C(A,A) \to HH_\idot(\C)
$$
compatible with the $HH^\hdot(\C)$-module structure on both sides.
\end{itemize}
The trace map also does not appear explicitly in \cite{kel}, but one
of its corollaries does. Namely, the trace map can be applied to the
identity map $\id_A$; the result is the {\em Chern class}
$$
\ch_A = \tr(\id_A) \in HH_0(\C).
$$
Keller's definition of the Chern class is slightly different, but
its main property is the same.
\begin{itemize}
\item for any short exact sequence $0 \to B \to C \to A \to 0$ in
  $\C$ we have $\ch_C = \ch_A + \ch_B$.
\end{itemize}
This is known as {\em devissage}. Other properties of Hochschild
homology are very prominent in \cite{kel}; among other things, they
insure that $HH_\idot(\C)$ is what it should be for particular
categories such as modules over an algebra or coherent sheaves on a
scheme. Of these properties, we will need the following two.
\begin{itemize}
\item If $\C_0 \subset \C$ is a thick (a.k.a. Serre) abelian
  subcategory, and $\C/\C_0$ is the quotient abelian category, then
  there exists a natural long exact sequence
$$
\begin{CD}
HH_\idot(\C_0) @>>> HH_\idot(\C) @>>> HH_\idot(\C/\C_0) @>>>
\end{CD}
$$
(the so-called {\em excision property}). If $\C$ is the category of
finitely generated modules over a Noetherian algebra $A$, then
$HH_\idot(\C) \cong HH_\idot(A) = \Tor_\idot^{A^{opp} \otimes
A}(A,A)$.
\end{itemize}
These properties are not used in the proof of the Gabber's Theorem
directly; their only importance is that they allow to compute
homology and cohomology of the relevant abelian category.

\section{Computations for sheaves.}\label{shv}

We will now sketch a direct -- and therefore, somewhat ugly --
construction of the Hochschild homology and cohomology theories in
the particular case needed for the Gabber's Theorem; it will satisfy
all the properties listed above. The category $\C$ in question is
the category $\Shv(X,Z)$ of sheaves on a smooth scheme $X$ supported
in a closed subset $Z \subset X$. We further assume that the scheme
$X$ is of finite type over a field $k$. We note that the case $Z =
X$, $\Shv(X,Z) = \Shv(X)$ is more-or-less classic by now, see
e.g. \cite{w}, \cite{S}. However, for the Gabber's Theorem, it is
essential to be able fix the support.

Consider the product $X \times_k X$, and let $\calo_\Delta$ be the
structure sheaf of the diagonal $\Delta \in X \times X$. Set
$HH^\hdot(\Shv(X)) = \Ext^\hdot(\calo_\Delta,\calo_\Delta)$.  We
note that these $\Ext$-groups can be computed by the Koszul
resolution. In particular, they remain the same if computed in the
category $\Shv(X \times X,\Delta)$ of coherent sheaves on $X \times
X$ supported on the diagonal $\Delta \subset X \times X$. Every
object $\F \in \Shv(X \times X,\Delta)$ defines a right-exact
functor $K(\F)$ from $\Shv(X)$ to itself by
$$
K(\F)(\E)=p_{2*}(p_1^*\E \otimes \F),
$$
where $p_{1,2}:X \times X \to X$ are the natural projections. The
correspondence $\F \mapsto K(\F)$ is an exact functor from $\Shv(X
\times X,\Delta)$ to $\Fun(\Shv(X),\Shv(X))$ which sends
$\calo_\Delta$ to the identity functor and induces therefore an
algegra map
$$
K:HH^\hdot(\Shv(X)) = \Ext^\hdot(\calo_\Delta,\calo_\Delta) \to
\Ext^\hdot(\Id_{\Shv(X)},\Id_{\Shv(X)}),
$$
which we take as our map $\ev$. If we restrict our attention to the
subcategory $\Shv(X,Z) \subset \Shv(X)$ supported in some closed $Z
\subset X$, then all of the above goes through literally, with one
change: one has to replace the product $X \times X$ with its formal
completion along $Z \times Z \subset X \times X$.

If $X$ is affine, then the algebra $HH^\hdot(\Shv(X))$ is easy to
compute explictly by using the Koszul resolution and the
local-to-global spectral sequence, which in these assumptions
degenerates: we have
$$
HH^\hdot(\Shv(X)) = H^0(X,\Lambda^\hdot\T(X)),
$$
where $\Lambda^\hdot$ is the exterior algebra, and $\T(X)$ is the
tangent sheaf of the scheme $X$ (since by assumption $X$ is regular,
$\T(X)$ is a vector bundle). Passing from $\Shv(X)$ to $\Shv(X,Z)$
amount to taking completion with respect to the ideal $\J_Z \subset
\calo_X$ defining $Z \subset X$. In particular, classes in
$HH^2(\Shv(X,Z))$ are just bivector fields on the formal completion
of $X$ along $Z \subset X$.

If $X = \Spec A$ is affine, then the deformations of the category
$\Shv(X)$ that we will need come from deformations of the algebra
$A$ (we have already seens the trivial deformation in
Example~\ref{exa}). Namely, by a first-order flat deformation of $A$
we will understand an associative algebra $A_h$ which is flat over
$k\langle h \rangle = k[h]/h^2$ and equipped with an isomorphism
$A_h/h \cong A$. Then for any such $A_h$, the category of
finitely-degenerated left $A_h$-modules obviously satisfies the
assumptions of Definition~\ref{cat.def}. One checks easily that for
any $a,b \in A$ lifted to elements $a',b' \in A_h$ we have
$$
a'b'-b'a'=h(\Theta_{A_h} \cntrct da \wedge ab)
$$
for some bivector field $\Theta_{A_h}$ on $A$, independent of $a$,
$b$, $a'$, $b'$. It is this bivector field $\Theta_{A_h} \in
HH^2(\Shv(X))$ that one associates to the deformation $A_h$. The
first piece of general nonsense that we accept without proof in this
paper is the following compatibility statement.

\begin{lemma}
The class $\Theta_{A_h}$ satisfies \eqref{def.comp}.\endproof
\end{lemma}

In other words, the (standard) definition of the deformation class
$\Theta_{A_h}$ agrees with the categorical definition. The proof is
rather traightforward but tedious: one replaces the functor categories
in the categorical definition by appropriate categories of
bimodules, and computes the class $\iota(t_{\C'})$ explcitly by
using the Koszul complex. We feel that in the context of this paper,
this is better left to the reader.

\bigskip

We now turn to the Hochschild homology. The groups
$HH_\idot(\Shv(X))$ are defined as
$$
HH_\idot(\Shv(X)) = H^\hdot(X,L^\hdot \delta^*\calo_\Delta),
$$
where $\delta:X \to X \times X$ is the diagonal embedding (in
algebraic notation, we have $L^\hdot \delta^*\calo_\Delta =
\Tor^\hdot_{X \times X}(\calo_\Delta,\calo_\Delta)$. We note that if
$X$ is not affine, $HH_\idot(\Shv(X))$ might be non-trivial both in
positive and in negative degrees. When we consider $\Shv(X,Z)
\subset \Shv(X)$, then the Excision property forces us to set
$$
HH_\idot(\Shv(X,Z)) = H^\hdot_Z(L^\hdot \delta^*\calo_\Delta),
$$
where $H^\hdot_Z$ denote cohomology with supports in $Z$. This might
be non-trivial both in positive and in negative degrees even when
$X$ is affine. For regular $X$, $L^\hdot \delta^*\calo_\Delta$ can
be computed explictly: the Hochschild-Kostant-Rosenberg Theorem
\cite{HKR} claims that
$$
L^i \delta^*\calo_\Delta \cong \Omega^i(X).
$$
Thus the homology sheaves of the complex $L^\hdot
\delta^*\calo_\Delta$ are sheaves $\Omega^i(X)$ of $i$-forms on $X$,
each placed in degree $-i$. If $X$ is affine, then the natural
$HH^\hdot(\Shv(X,Z))$-structure on $HH_\idot(\Shv(X,Z))$ is given by
the usual contration operation between a form and a polyvector
field. We note, and this is important, that the complex $L^\hdot
\delta^*\calo_\Delta$ on $X$ is Serre self-dual.

The most difficult thing to define explicitly is the trace map. One
approach is the following. Let $\E \in \Shv(X)$ be a coherent sheaf
on $X$. Then by adjunction, we have a natural class
$$
a_\E \in \Ext^0(L^\hdot \delta^*\delta_*\E,\E) \cong \Ext^0(\E
\otimes L^\hdot \delta^*\calo_\Delta,\E) \cong \Ext^0(\E,\E \otimes
L^\hdot \delta^*\calo_\Delta),
$$
where the first isomorphism follows from the projection formula, and
the second one uses the fact that $L^\hdot \delta^*\calo_\Delta$ is
Serre self-dual. This class $a_\E$, while quite tautological, is not
trivial; in characteristic $0$, it incorporates the Atiyah class $A
\in \Ext^1(\E,\E \otimes \Omega^1(X))$ together with its exterior
powers $\Lambda^k(A) \in \Ext^k(\E,\E \otimes \Omega^k(X))$ (see
\cite{mar} for a detailed exposition). We note that while one
usually considers the Atiyah class for vector bundles, there is no
need to do so: the definition works just as well for every coherent
sheaf. By transposition, for every coherent sheaf $\E$ we obtain a
class
$$
a_\E^t \in \Ext^0(\rhom^\hdot(\E,\E),L^\hdot \delta^*\calo_\Delta),
$$
where $\rhom^\hdot(\E,\E)$ is the complex which represents the local
$\Ext$-groups from $\E$ to itself. Now, assume that $\E$ is in fact
supported in $Z \subset X$. Then so is the complex
$\rhom(\E,\E)$. Therefore, by taking cohomology with supports in $Z
\subset X$, the class $a^t$ gives a map
$$
H^\hdot_Z(\rhom^\hdot(\E,\E)) \cong \RHom^\hdot(\E,\E) \to
HH_\idot(\Shv(X,Z)).
$$
This is our trace map $\tr_{\E}$. We note that by its very
construction, it is compatible with the natural
$HH^\hdot(\Shv(X,Z))$-module structures on both sides. Checking
devissage for the associated Chern character $\ch_{\E} =
\tr_{\E}(\id_{\E})$ is also straighforward. Indeed, for any short
exact sequence
$$
\begin{CD}
0 @>>> \E_1 @>>> \E @>>> \E_2 @>>> 0
\end{CD}
$$
of sheaves on $X$, we can consider $\E$ as a filtered sheaf (with a
two-step filtration $F_0\E=\E_1$, $F_1\E=\E$). Then the algebra
$\rhom^\hdot_F(\E,\E)$ of filtered local $\Ext$-groups is
well-defined, and we have the commutative diagram
$$
\begin{CD}
\rhom^\hdot_F(\E,\E) @>>> \rhom^\hdot(\E,\E)\\
@VVV  @V{a^t_{\E}}VV\\
\rhom^\hdot(\E_1) \oplus \rhom^\hdot(\E_2) @>{a^t_{\E_1} \oplus
  a^t_{\E_2}}>> L^\hdot \delta^*\calo_\Delta
\end{CD}.
$$
It remains to notice that $\Id_\E$ comes from the class $\Id^F_\E
\in \Ext^0_F(\E,\E)$ which projects to the class $\Id_{\E_1} \oplus
\Id_{\E_2} \in \Ext^0(\E_1,\E_1)\oplus\Ext^0(\E_1,\E_1)$.

The second general-nonsense compatibility result that we prefer to
accept without any indication of proof concern an explicit
computation of the Chern character in one simple case. Assume that
the regular scheme $X$ is affine, and moreover, assume that $Z
\subset X$ is also regular, of codimension $l$. Then the local
cohomology module $H^p_Z(X,\calo_X)$ is trivial for $p \neq l$, and
$H^l_Z(X,\calo_X)$ is well-known: this is the so-called (right)
$\delta$-function $D$-module on $X$ supported in $Z$. It has an
increasing filtration whose associated graded quotient is isomorphic
to
$$
\gr^\hdot H^l_Z(\calo_X) \cong H^0(Z,\omega_{X,Z}^{-1} \otimes
S^\hdot\N_{X,Z}),
$$
where $\N_{X,Z}$ is the normal bundle to $Z$ in $X$, $S^\hdot$
stands for symmetric power, and $\omega_{X,Z} =
\Lambda^l(\N^*_{X,Z})$ is the top exterior power of the conormal
bundle $\N^*_{X,Z}$. In particular, we have a canonical embedding
$H^0(Z,\omega^{-1}_{X,Z}) \to H^l_Z(X,\calo_X)$. This gives
canonical maps
$$
H^0(Z,\omega^{-1}_{X,Z} \otimes \Omega^p(X)|_Z) \to
HH_{p-l}(\Shv(X,Z))
$$
for any $p$, $0 \leq p \leq \dim X$. But for $p=l$, the natural map
$\omega_{X,Z} \cong \Lambda^l\N^*_{X,Z} \to \Omega^l(X)|_Z$ gives a
tautological class $\tau_{X,Z}$ in $H^0(Z,\omega^{-1}_{X,Z} \otimes
\Omega^l(X)|_Z)$ and consequently, in $HH_0(\Shv(X,Z))$.

\begin{lemma}\label{tau}
In the assumptions above, the Chern character
$$
\ch_{\calo_Z} \in HH_0(\Shv(X,Z))
$$
coincides with the tautological class $\tau_{X,Z}$.\endproof
\end{lemma}

\begin{remark}\label{tr.rem}
The real problem in this compatibility Lemma is to use a convenient
construction of the embedding $H^0(Z,\omega^{-1}_{X,Z}) \to
H^l_Z(X,\calo_X)$. In practice, this embedding is usually obtained
by using residues, as in \cite{hart}. For purely formal reasons such
as Excision, any well-developed Hochschild homology formalism would
contain this as an integral part, in the form of the Tate residue
and its higher-dimensional generalizations. Then Lemma~\ref{tau}
would essentially be the definition of the class $\tau_{X,Z}$; the
only thing to prove would be $\tau_{X,Z} \neq 0$, and even this
should be an easy exercize. Unfortunately, Tate residue is not part
of the foundational paper \cite{kel} (this is exactly one of those
applied results that are likely to be omitted at the first stages of
development of the general theory).
\end{remark}

\section{Gabber's Theorem.}\label{pf}

Having spent several pages, for better or for worse, on describing
the general formalism of Hochschild Homology and Cohomology, we can
now show how it helps to prove Gabber's Theorem. The proof is really
short.

By the preliminary reductions done in \cite{G}, the theorem amounts
to the following purely algebraic fact. Assume given a commutative
algebra $A$ over a field $k$ of characteristic $0$, a prime ideal
$\m \subset A$, and a finitely generated $A$-module $M$ annihilated
by a power of $\m$. Assume that $A$ and $A/\m$ are Noetherian and
regular. Assume in addition that we are given a first-order
one-parameter deformation $A_h$ of the algebra $A$, and assume that
there exists a left $A_h$-module $M_h$ which is flat over $k[h]/h^2$
and satisfies $M_h/h \cong M$. Define a bracket operation $\{-,-\}$
on $A$ by
$$
a'b'-b'a' = h\{a,b\}
$$
for any $a,b \in A$, where $a', b' \in A_h$ are arbitrary liftings
of the elements $a,b$. Since $A_h/h\cong A$ is commutative, this
does not depend on the choice of liftings; if $A_h$ comes from a
full one-parameter non-commutative deformation of $A$, then
$\{-,-\}$ is the associated Poisson bracket.

\begin{theorem}
In the assumptions above, the ideal $\m \subset A$ is involutive,
$\{\m,\m\} \subset \m$.
\end{theorem}

\proof{} Let $\Theta_{A_h}=\Theta$ be the bivector field on $\Spec
A$ associated to the deformation $A_h$, so that $\{a,b\}=\Theta
\cntrct (da \wedge db)$. Then obviously $\{\m,\m\}\subset\m$ if and
only if the $\eta(\Theta)=0$, where $\eta:\Lambda^2\T(A) \to
\Lambda^2\N_{A,\m}$ is the projection onto the second exterior power
of the normal bundle $\N_{A,\m}$. This is equivalent to
\begin{equation}\label{cnt}
\Theta \cntrct \omega = 0,
\end{equation}
where $\omega \in \Omega^\hdot(A)/\m$ is the determinant of the
conormal bundle $\N^*_{A,\m} \cong \m/\m^2$. Take now $X = \Spec A$,
and let $Z \subset X$ be defined by the ideal $\m$. Then by
Lemma~\ref{tau}, \eqref{cnt} is equivalent to
$$
\Theta \cdot \ch_{A/\m} = 0 \in HH_{-2}(\Shv(X,Z)).
$$
To prove this, we note that $HH_{-2}(\Shv(X,Z))$ is a flat module
over $A/\m$, so that we can localize $A$ at $\m$ and assume that $\m
\subset A$ is a maximal ideal. Then on one hand, $\Theta \cdot \ch_M
= \Theta \cdot \tr_M(\id_M) = \tr_M(\Theta \cdot \id_M) =
\tr_M(\ev_M)$, which is $0$ by Lemma~\ref{no.ev}, and on the other
hand, by devissage
$$
\ch_M = n\ch_{A/\m}
$$
for some positive integer $n$.
\endproof

In a nutshell, the essential part of this proof is this: we want to
reduce the claim for general $M$ to the (simple) particular case
when $M$ is annihilated by $\m$, not some power $\m^l$. We do this
by reducing the claim to statement about a Hochschild homology
class, and then using devissage to replace $M$ with the associated
graded quotient $\gr M$ with respect to the $\m$-adic filtration.
Just as \cite{GS} and \cite{G}, our proof is essentially a ``trace
argument''; trace appears in the form of the Chern class map.

We note that the argument is extremely crude -- localizing at $\m$,
the generic point of $Z \subset X$, we lose a lot of
information. When one has a well-developed Hochschild homology
formalism at one's disposal, the claim could be considerably
strengthened to obtain information about the Chern class of the
module $M$ (and the associated sheaf, if one no longer assumes that
$X$ is affine). In particular, it should be possible to obtain a
general $D$-module version of the obvious fact that the vector
bundle underlying a local system has trivial Chern classes. All
this, however, remains a topic for future research.

One point in which our proof is seriously deficient in comparison to
\cite{G} is our assumption that $X$ is regular. In $D$-module
applications, this is a given; in fact, the $D$-module theory on
singular varieties is developed in a different way, and Gabber's
Theorem for singular $A$ would probably be useless for
this. However, the Theorem itself holds in full generality -- no
assumptions of regularity are imposed in \cite{G}. In our argument,
they are not strictly needed, either. However, one has a feeling
that a reasonable theory of Hochschild Homology can only be
developed for categories of finite homological dimension. Thus for
general $A$, our general nonsense statements are probably not true
without serious modifications, and anything built on them becomes
suspect.

\bigskip

\noindent
{\sc Steklov Math Institute\\
Moscow, USSR}

\bigskip

\noindent
{\em E-mail address\/}: {\tt kaledin@mccme.ru}

\end{document}